\documentclass[12pt,a4paper]{article}
\usepackage{color}
\usepackage{setspace}
\usepackage{fancyvrb}
\usepackage{amssymb}
\usepackage[leqno]{amsmath}
\usepackage{graphicx}
\DeclareFontFamily{OT1}{codefont}{}
\DeclareFontShape{OT1}{codefont}{m}{n}{ <-> cmtt12 }{}

\title{A coefficient related to splay-to-root
traversal, correct to thousands of decimal places}

\author{Colm \'O D\'unlaing\thanks{e-mail: odunlain@maths.tcd.ie.
Mathematics department website: http://www.maths.tcd.ie.}\\
{\em Mathematics, Trinity College, Dublin 2, Ireland}}
\begin{document}
\renewcommand{\theequation}{\arabic{equation}}
\newtheorem{theorem}{Theorem}
\newtheorem{definition}[theorem]{Definition}
\newtheorem{lemma}[theorem]{Lemma}
\newtheorem{openproblem}[theorem]{Open Problem}
\newtheorem{remark}[theorem]{Remark}
\newtheorem{proposition}[theorem]{Proposition}
\newtheorem{corollary}[theorem]{Corollary}
\newtheorem{assumption}[theorem]{Assumption}
\newtheorem{specification}[theorem]{Specification}
\newtheorem{notation}[theorem]{Notation}
\newtheorem{hypothesis}[theorem]{Working Hypothesis}
\newcommand{\cost}{\text{\rm cost}}
\newcommand{\height}{\text{\rm height}}
\newcommand{\spine}{\text{\rm spine}}
\newcommand{\fetch}{\text{\rm fetch}}
\newcommand{\irp}{\text{\rm irp}}
\newcommand{\rp}{{\text{\rm rp}}}
\newcommand{\isl}{\text{\rm isl}}
\newcommand{\tsl}{\text{\rm tsl}}
\newcommand{\spl}{\text{\rm sl}}
\newcommand{\hp}{{\rm hp}}
\newcommand{\xp}{\text{\rm xp}}
\newcommand{\be}{\begin{equation*}}
\newcommand{\ee}{\end{equation*}}
\newcommand{\bl}{\begin{lemma}}
\newcommand{\el}{\end{lemma}}
\newcommand{\ol}{\overline}
\newcommand{\ul}{\underline}
\newcommand{\bs}{\backslash}
\newcommand{\fs}{/}
\newcommand{\n}{\emptyset}
\newcommand{\offline}{\text{\rm offline}}
\newcommand{\hb}{\hfil\break}
\newcommand{\caret}{\char`\^}
\newcommand{\stk}{{\text{\rm stk}}}
\newcommand{\pst}{{\text{\rm pst}}}
\newcommand{\pers}{{\text{\rm persistence}}}
\newcommand{\rmp}{{\text{\rm rmp}}}
\newcommand{\leqin}{{~\leq_{\text{\rm in}}~}}
\newcommand{\IR}{{\mathbb R}}
\def\qedrule{\vrule height0.75em width0.5em depth0.25em}
\def\makeqedrule{\nobreak \ifvmode \relax \else
      \ifdim\lastskip<1.5em \hskip-\lastskip
            \hskip1.5em plus0em minus0.5em \fi \nobreak
      \qedrule\fi}
\let\qed=\makeqedrule
\newcommand{\takeanumber}{\refstepcounter{theorem}%
\def\nextnumber{\thetheorem}}
\newcommand{\numpara}{\refstepcounter{theorem}%
\medbreak\noindent{\bf (\thetheorem)\ \ \ }}

\DeclareGraphicsExtensions{.eps}

\maketitle

\begin{abstract}
This paper takes another look at the cost of traversing
a binary tree using repeated splay-to-root.  This was
shown to cost $O(n)$ (in rotations) by Tarjan [\ref{tarjan85}]
and later, in different ways, by others [\ref{elmasry}].

It would be interesting to know the minimal possible coefficient
implied by the $O(n)$ cost; call this coefficient $\beta$.
In this paper we define a related coefficient $\alpha$
describing the cost of splay-to-root traversal on maximal (i.e.,
complete) binary trees, and show that $\beta \geq 2 + \alpha$.
We give the first 3009 digits of $\alpha$, including the decimal
point, and show that every digit is correct.

We make two conjectures: first, that $\beta = 2 + \alpha$, and
second, that $\alpha$ is irrational.
\end{abstract}

\section{Introduction}

In this paper, `tree' means `binary tree.'  The size of a 
tree $T$, the  number of nodes in the tree, is denoted $|T|$ and $\emptyset$
denotes the empty tree, of size zero.  The splay operations, `zig' (1 rotation),
`zigzig' and `zigzag' (2 rotations) were introduced in 
[\ref{sleator-tarjan}], and shown to lead to optimal amortised
costs for several operations.

This paper takes another look at the cost of traversing a tree using
repeated splay-to-root.
By $\cost(T)$ is meant the number of rotations in
a complete splay-to-root traversal of a tree $T$.
Tarjan [\ref{tarjan85}] showed that
$\cost(T)$ is $O(|T|)$.  Elmasry [\ref{elmasry}] gave a very
elegant derivation of this result,  giving a concrete upper bound of
$4.5 |T|$, but it appears that this estimate counts {\em links}
(splay operations)
rather than rotations, and the implicit bound, counting rotations,
would be $8|T|$.

An interesting problem is to determine the exact value of $\beta$, where
\be
\beta = \inf \left \{ \beta' \in \IR:~ ( \forall \epsilon > 0 )
(\exists N) (\forall T \not= \emptyset) \left ( |T| \geq N \implies \left (
\frac{\cost(T)}{|T|} < \beta' + \epsilon \right )\right )\right \}.
\ee

This infimum $\beta$ exists [\ref{tarjan85},\ref{elmasry}] and $\beta \leq 8$
[\ref{elmasry}].

We were able to calculate not $\beta$,
but a related constant  $\alpha$, correct to several
thousand decimal places, and provably so. Here
\be
\alpha = \lim_{h\to \infty} \frac{\cost(M_h)}{|M_h|},
\ee
where $M_h$ is the maximal tree of height $h$
(Definition \ref{def: maximal tree}) and the
cost is the number of rotations in the splay-to-root
traversal of $M_h$.  Actually, $\alpha$ is studied
as the limit of a slightly different sequence
\be
\alpha = \lim_h \frac{\cost(M_h)+1}{|M_h|+1}
= \lim_h \frac{\tsl(M_h)}{2^{h+1}}
\ee
(Definition \ref{def: alpha second version}; $\tsl$
is `total spine  length.').

What makes this interesting is that 
$\alpha$ is almost certainly irrational,
though proofs of irrationality are notoriously
difficult [\ref{kristensen}].

We show also that $\beta \geq 2+\alpha$.

Our analysis is different from the traditional $O(\ldots)$
because we aim at an accurate determination of $\alpha$.
The analysis mostly involves the definition, convergence,
and estimation of various series.

We make two conjectures: (i) that $\beta = 2+\alpha$ and
(ii) that $\alpha$ is irrational.

\subsection{Computing notes}

The number $\alpha_{10049}$ was originally computed in 60
hours using lazy evaluation of `root persistence' of
$M_h^{[1]}$, $0\leq h \leq 10498$,
all to be defined later.
There is nothing magical about 10049: the program was
stopped after sixty hours, that is all.

Actually, the root persistence
of the relevant $M_h^{[2]}$, and then $M_h^{[1]}$, can be computed in
22 seconds.  There would be no difficulty in producing
$\alpha_N$ for much larger $N$, possibly in the millions.

It seems unlikely that this would shed much light on the
question of periodicity of $\alpha$, though it might be
of interest to study the statistical distribution of digits.

\section{The constant $\alpha$}
Here is $\alpha$ to 3007 decimal places, and every
digit is correct.

\begin{Verbatim}
2.41464532311342664135721059929950736447077229680868005373357354525826
2740624105270973536945896664894612383387133140143583762512616478135561
4286064288849821865109967561998433943559594398275126178907483349665025
5080091545859089478186004505638394411483898674045412979196509367417534
6878187831010615663204826748679513555160671568781383651194277803644514
9238821898220035037296172434065926647433143400078616634513025688467682
6171027174616776487580847038386785938892724884200227357265301260112929
0242872683045128992893355401230833822607056585290386388421722662334074
0116380336801125542664867713593347299299395197919439890378605878338449
0684574027066720810192492640915295144745925864171925599305806057575314
2508105256778113408318294187577539205456236837958433846424262409916448
9697233372775761982930647206283833774577192257755235647190176672592133
2949875562869242099200719629819461975554171376684916338003857643161270
0747640412245784792115264723755558116173763305314958459449252363780334
3121387447870909378536449539445023339622338606613085762454506694594937
1377070000570876683583128296299442896491275160724353720907079555723091
7471898153816752661145112486336096680922635899639696366033792702634139
2837615929543809383793534957141012590712487704983953605803254154040608
8106227758784127307004889163446215318322212040230203240483964362931477
7176498520706885768943453560729740034369991068278397299835094676979973
0685469380069957582763661847207759482022362106606817207185249175670881
9473275669722948259829248370547849180846820987183968707971251129981407
9773527676551260720609151582081068413751466710205527377054258149097152
9289561262502780864046114682107387796782939836021136317806326521145487
6760227391835198547411927129861199542826677693150395041342412196574665
4473060865865080514373454048919553626044024615908890117833068246524507
0089664144159195043023934386059335422404227281126486873705576624384602
9528750910516010997873931055591565295205182888286563627909743235573390
9689133141395995283844816911095381106026450878554652596426441091822593
2390806913732320816254273685041076786682097197712812201979942001490611
7626168639781158815517923973646158921003560714955895298284159796312644
0396095898620587149904372166146231321103503671310912936678042493588535
2505353868931579424411182910431741623568322771304576555905031402969234
4660957154696322464495584555868741912213583503440498876463365928113074
9233658183765301672791687233629880536085398659358936376703856335909377
6948499772272644410353046707274139060905827949996769797881256122388937
2692185590090760815441293071672841610891695407335595431360249781445226
4120383844036525829611612062383372947896234386973096119274358253660018
7150107446002431807957835268023887676251646716118886314819424666594719
2923138902667033619091930830078377030802184358525745880030925494397385
5546393947770581175494013381395303045868156058494995588030890182315469
2099632400578717275003880135554963167113109913616699600867929624278862
052537884913130875294295119205291521849749337282042856555796587734387
\end{Verbatim}

\begin{lemma}
The above is correct in all 3007 decimal places.
\end{lemma}

{\bf Proof.}
We know that 
\be
\alpha \leq \alpha_{10049} + 8\times \frac{10053^4}{2^{10049}}
\ee
(Corollary \ref{cor: main inequality}). The tolerance, that
is, the number
$8 \times \frac{10053^4}{2^{10049}}$, was computed exactly
(10050 decimal places), then truncated to 3150 characters.
The last 3 lines (210 digits) of the tolerance are:

\begin{Verbatim}
0000000000000000000000000000000000000000000000000000000000000000000000
9094028331545896961083851103178971716996186990956575118754540507007346
6666935419158412707593514820820039606856557926859110469901790668773529
\end{Verbatim}

Also, $\alpha_{10049}$ was calculated exactly.  Truncated, the last 3 lines are
\begin{Verbatim}
0525378849131308752942951192052915218497493372820428565557965877343873
                                                                    *
9248328553474208432836925319209482543964256661337568102900090351978867
9254269613807704859718778867190874549841852135271086734458892430187008
\end{Verbatim}
This is a lower bound.
The sum of $\alpha_{10049}$ and the tolerance was calculated exactly,
and again truncated, to give an upper bound whose last 3 lines are
\begin{Verbatim}
0525378849131308752942951192052915218497493372820428565557965877343874
                                                                    *
0157731386628798128945310429527379715663875360433225614775544402679602
5920963155723546130478130349272878510527507927956997781449071497064361
\end{Verbatim}
It is evident that the lower and upper bound agree up the point
marked: hence the result.\qed

\subsection{Is $\alpha$ irrational?}

We believe so because it {\em looks} irrational.
If it is rational, it is periodic.  
The methods of Knuth, Morris, and Pratt were used in a na\"{\i}ve
way to look for self-overlaps, from right to left. This
showed a maximum self-overlap length of $3$, which means
that there are no periods.

\section{Splay-to-root traversal}
\subsection{Various definitions}
The coefficient $\alpha$
describes the traversal cost of {\em maximal} trees.
We use the word `maximal' in preference to `complete'
because the latter is ambiguous.

\begin{definition}
\label{def: maximal tree}
The {\em height} of a nonempty tree $T$ is the number of links
(not nodes) in the longest path from the root to a leaf.
The empty tree has height $-1$.

A tree $T$ of height $h$ is {\em maximal} if it has the
maximal possible number of nodes, $2^{h+1}-1$, for trees
of height $h$.  For each $h\geq -1$,
there is exactly one maximal tree of
height $h$.  We denote by
\be
M_h
\ee
the maximal tree of height $h$ ($M_{-1}$ can be identified with the
empty tree).
\end{definition}

Then $\alpha$ can be described as follows, though we shall
work with an
equivalent version, introduced later.

\begin{definition}
\label{def: alpha first version}
The coefficient $\alpha$: first version.
\be
\alpha = \lim_{h\to \infty} \frac{\cost(M_h)}{|M_h|}
\ee
\end{definition}

It is easier to work with `fetch and discard,'
introduced below, instead of `splay to root.'  Both
procedures have the same cost.

\subsection{Spine lengths}

In traversing a tree $T$ by repeated splay-to-root, at
the first step the leftmost node in $T$ is brought to the
root, or `fetched.'  At that time the leftmost branch from
the root is called the {\em spine.} The {\em length} of the
spine is the number of {\em nodes} on the spine.  The
cost (in rotations) of the first fetch is the spine length minus 1.

In all subsequent steps, the spine is the leftmost branch from
the right child of the root, and since that child has depth 1,
the  cost of fetching equals the length of the spine.

\begin{corollary}
Given a nonempty tree $T$, $\cost(T)$,
the cost in rotations of traversing
$T$ by repeated splay-to-root, is the total spine length minus 1.\qed
\end{corollary}

\subsection{Fetch and discard}
\label{subsect: fetch and discard}

\begin{figure}
\centerline{
\includegraphics[width=4in]{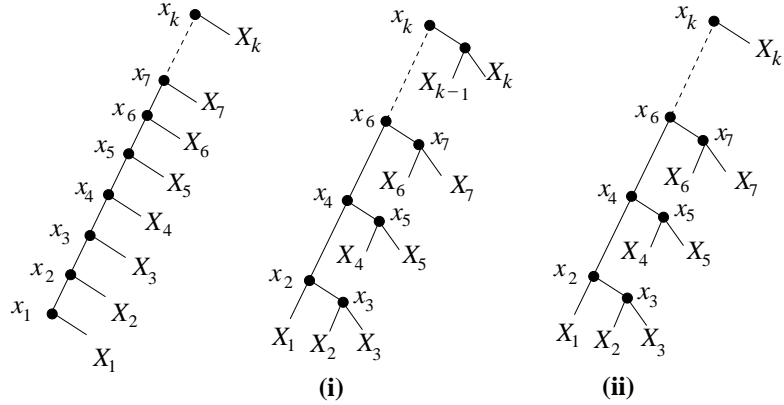}
}
\caption{One step of fetch-and-discard; $k$ is the spine length.
(i) $k$ odd; (ii) $k$ even.}
\label{fig: f-a-d}
\end{figure}

This is a modification of splay-to-root where
every time a node is brought to the root by fetching, it
is deleted from the tree.
If $T'$ is the result of $k$ steps of splay-to-root,
and $T''$ is the result of $k$ steps of fetch-and-discard,
where $k > 0$, then $T''$ is isomorphic to the right subtree of the root
in $T'$.

Therefore the {\em spine} in fetch-and-discard, i.e., the
leftmost branch containing the next node to be fetched, is always the
leftmost branch from the root.

\subsection{Fetch and discard, in detail}

The effect of fetching on
a tree $T$ is illustrated in Figure \ref{fig: f-a-d}.
The spine is labelled $x_1,\ldots, x_k$, from bottom to top.
Also, $X_j$ is the right subtree  of $x_j$, $1 \leq j \leq k$.
The effect of fetching the leftmost node
$x_1$ by splay operations is as follows.

\begin{itemize}
\item
If $k=1$ then $T$ is replaced by $X_1$ (i.e., the
root of $T$ is the root of $X_1$) and the operation is
finished. Otherwise, $k>1$.
\item
For $j = 2,4,\ldots$, if $j < k$,
$x_{j+1}$ is `pushed off the spine' in the sense that
$x_j$ remains on the spine, and $x_{j+1}$ becomes the
right child of $x_j$, $X_j$ becomes the left subtree
of $x_{j+1}$, and $X_{j+1}$  continues as the right subtree
of $x_{j+1}$.
\item
If $k$ is even, then $x_k$ remains on top of the spine
and its right subtree continues as $X_k$.  If $k$ is odd
(and $k>1$), then $x_k$ has been made the right child of
$x_{k-1}$.
\item
$X_1$ becomes the left subtree of $x_2$.
\item
$x_1$ is discarded.
\end{itemize}

\begin{definition}
{\rm (a)} Given $0 \leq k \leq |T|$,
\be
\fetch(k,T)
\ee is the tree which results
when fetch and discard is applied to $T$ $k$ times.

\hb
\verb\      \{\rm (b)} $|\spine(T)|$ is the spine length,
the number of nodes in $\spine(T)$, and

\begin{gather*}
\text{\rm (c)}\quad \tsl(T) = 
\sum_{k\leq |T|} |\spine(\fetch(k,T))|.
\end{gather*}
(In other words, $\tsl(T)$ is the total spine length in
fetch-and-discard.)
\end{definition}

\begin{corollary}
\label{cor: cost = tsl-1}
If $T$ is a nonempty tree, then $\cost(T)$, the
cost in rotations of splay-to-root traversal of $T$, satisfies
$\cost(T) = \tsl(T)-1$.\qed
\end{corollary}

\begin{figure}
\centerline{
\includegraphics[width=5in]{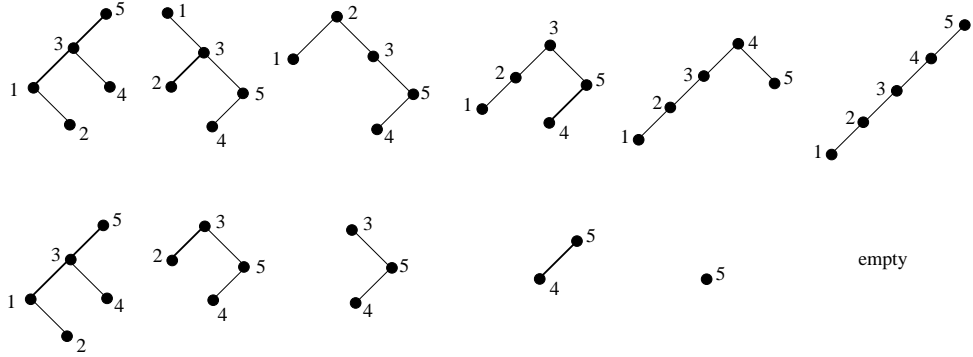}}
\caption{Comparing splay-to-root with fetch-and-discard}
\label{fig: discard}
\end{figure}

Figure \ref{fig: discard} shows how traversal by splay to root
and by fetch and discard have the same cost, if one counts
total spine length.

It is convenient to use $\tsl(T)$, rather than rotation
count, as the cost of
traversing a tree by fetch and discard.
The difference, $1$ if $T$ is nonempty, is negligible.

From Definition \ref{def: alpha first version}, and
Corollary \ref{cor: cost = tsl-1}, and since $|M_h| = 2^{h+1}-1$,
\be
\alpha = \lim_{h\to \infty} \frac{\tsl(M_h)-1}{2^{h+1}-1} .
\ee

In fact, we shall use the following equivalent definition
of $\alpha$.  It is easier to handle. Both the numerator
and denominator have been increased by $1$.  We have
yet to show that the limit exists.

\begin{definition}
\label{def: alpha second version}
The coefficient $\alpha$: second, and equivalent, version.
\begin{gather*}
\alpha = \lim_{h \to \infty} \alpha_h,\quad\text{where}\\
\alpha_h = \frac{\tsl(M_h)}{2^{h+1}}
\end{gather*}
\end{definition}

\section{$\beta \geq 2+\alpha$}
\begin{figure}
\centerline{\includegraphics[height=2in]{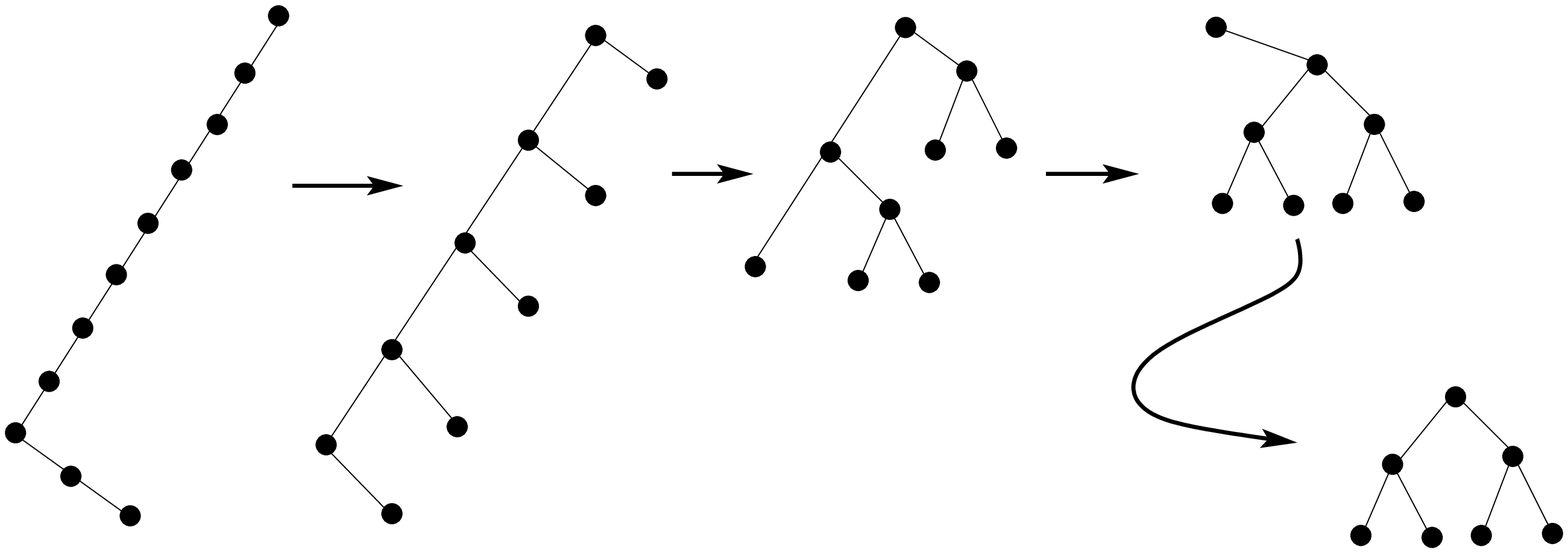}}
\caption{$2^r+r$ nodes reduce to $M_{r-1}$.}
\label{fig: beta2alpha}
\end{figure}

\begin{theorem}
\label{thm: beta2alpha}
$\beta \geq 2+\alpha$
\end{theorem}

{\bf Proof.}  Let $T$ be a tree of $2^r+r$ nodes, whose
leftmost node is on a spine of $2^r+1$ nodes and from
which there is a rightmost branch of $r$ nodes.

Now, $r$ fetches will reduce spine to a single node
whose right subtree is $M_{r-1}$.  One more
fetch will reduce the tree to $M_{r-1}$.

The total spine length in these $r+1$ fetches is
\be
2^r+1 + 2^{r-1}+1 + \ldots + 2 + 1 + 1
= r + 2^{r+1} - 1.
\ee
Therefore 
\be
\tsl(T) = 2^{r+1} + r-1 + \tsl(M_{r-1})
\ee
Divide by $2^r$:
\be
\left ( \frac{\tsl(T)}{2^r+r}\right )
\left (  1 + \frac{r}{2^r} \right ) = 2 + \frac{r-1}{2^r} + \alpha_{r-1}
\ee
For large $r$, the left-hand side
is close to $\tsl(T)/|T|$ and the right-hand side is close to $2+\alpha_{r-1}$.
See Figure \ref{fig: beta2alpha}.\qed

\section{Upper segments}
\begin{figure}
\centerline{\includegraphics[height=1in]{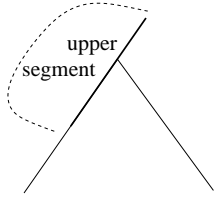}}
\caption{Upper segment marked by the heavy line.}
\label{fig: upper-seg}
\end{figure}

The following will be used in Lemma \ref{lem: continuation}.

\begin{definition}
\label{def: upper segment}
An {\em upper segment} of a tree $T$ is a subset $E$ of
$\spine(T)$ with the property that if $u\in E$ and $u$ is
not the root then the parent of $u$ is also  in $E$.
See Figure \ref{fig: upper-seg}.
\end{definition}

\begin{figure}
\centerline{\includegraphics[height=1in]{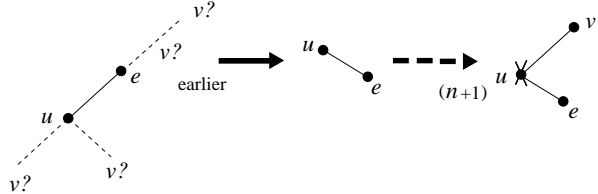}}
\caption{Lemma \ref{lem: preserve upper segments}, case (i)
illustrated.  The argument is that $v$ was above $e$ at
an earlier step.}
\label{fig: upperseg-i}
\end{figure}

\begin{lemma}
\label{lem: preserve upper segments}
Suppose that $E$ is an upper segment of a tree $A$
subject to $n$ fetch-and-discard steps.
Then

\centerline{
{\rm (*)}  for every subtree $T$ of $\fetch(n,A)$,
$E\cap \spine(T)$
is an upper segment of $T$.
}
\end{lemma}

{\bf Proof.} By induction on the number $n$ of fetch operations.

Suppose that the first $n$ operations preserve the condition (*)
but the $n+1$-st does not.

Let $B=\fetch(n,A)$.
There are two possibilities.

(i) The first possibility is that a spine node $u$ in $B$ has parent $v$ and right child $e$
where $e \in E$ and $v\notin E$, and $v$ is pushed off the spine.
Thus $v$, not in $E$, acquires a left child which is in $E$.
Since $e$ was a right child of $u$, there must have been
an earlier step when $e$ was parent of $u$ on the spine and became its
right child.

We claim that at this earlier step, $v$ was above $e$ on the spine.
See Figure \ref{fig: upperseg-i}.

In support of this claim,
(a) $v$ cannot have been in the left subtree at $u$,
since it follows $u$ in inorder;  (b) $v$ cannot
have been in the right subtree at $u$ since it would not
reach the spine until $u$ was fetched;  (c) $v$
cannot have been in the right subtree of any spine node
above $u$, since again it would not reach the spine until
after $u$ was fetched.  So $v$ is above $u$ on the spine,
and since $e$ was the parent of $u$ at this time, $v$
is above $e$ on the spine, as claimed.
Then (*) was violated at an earlier step.

(ii) The other possibility is that
the node $x$ being fetched from $B$ has right subtree
$T$, the root $e$ of $T$ is in $E$, and the parent $v$ of
$x$ (on the spine) is not in $E$.  Again, $e$ must have
been pushed off the spine by $x$ at an earler step.

(a) $v$ cannot have been in the left subtree at $x$ since
it would have been fetched before $x$; (b) $v$ cannot have
been in the right subtree at $x$, since it would remain
off the spine until $x$ is fetched; and (c) $v$ cannot
be in the right subtree of any other spine node since
it would remain off the spine until $x$ is fetched.  Therefore,
$v$ was above $x$ when $x$ pushed $e$ off the spine,
so $e$ was below $v$ on the spine, $e \in E$, $v\notin E$,
and (*) was violated at an earlier step.\qed

\section{Extending and combining trees}
The two results in this section
applicable to estimating $\tsl(M_h)$ and hence $\alpha$, are
Corollary \ref{cor: tsl mh+1 rp mh1} and
Corollary \ref{cor: 2nd diff}.

\subsection{Extensions}
\begin{definition}
Given two trees $A,B$ and a node $x$ not in $A$ or $B$,
$A^xB$ is the tree with root $x$ and left and right subtrees
$A$ and $B$.  Also, $A^x = A^x \emptyset$: that is,
$A^x$ has root $x$, and $x$ has left subtree $A$
and empty right subtree: $x$ is the rightmost node in inorder.
\end{definition}

\begin{figure}
\centerline{
\includegraphics[height=2in]{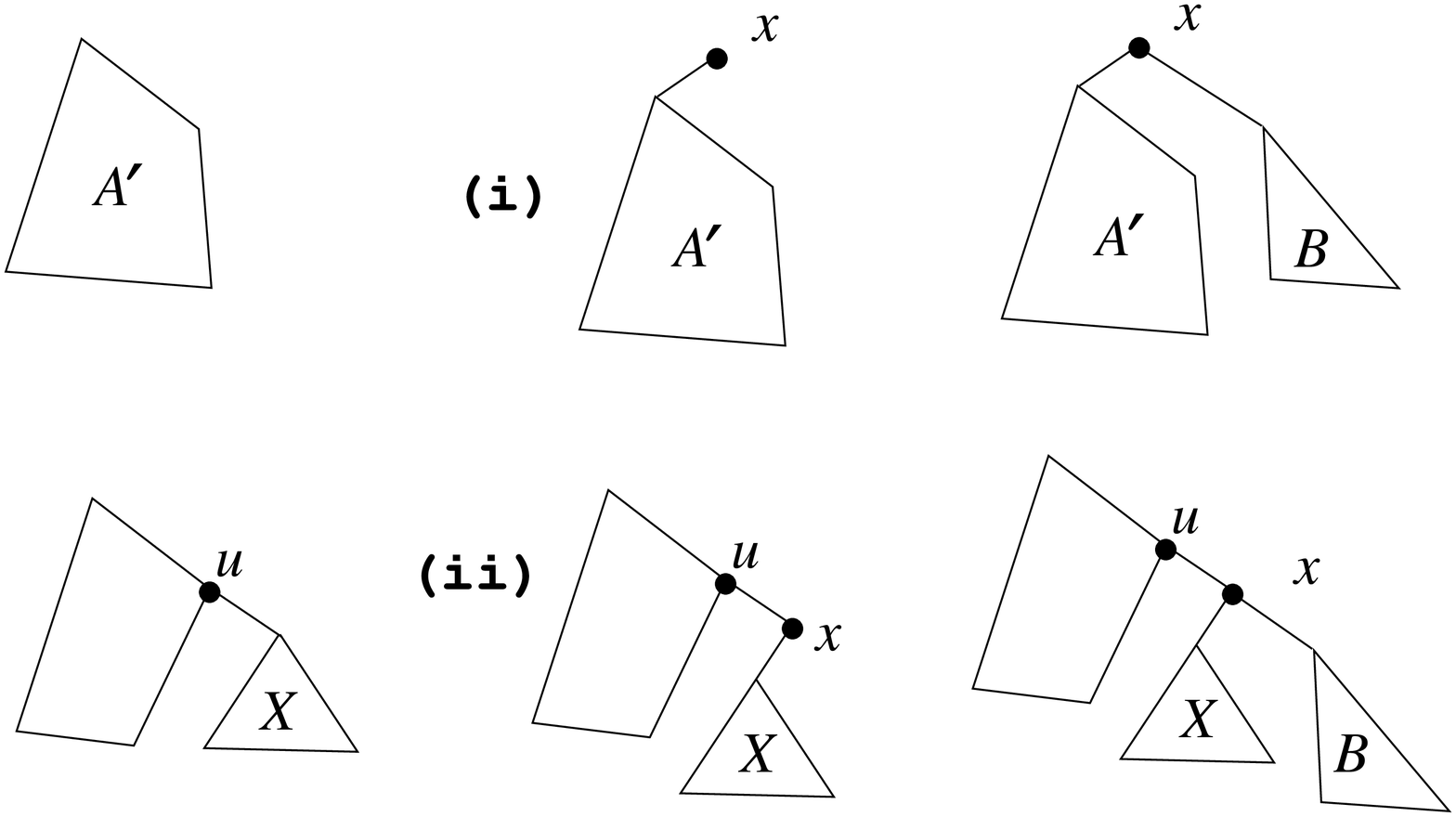}
}
\caption{Illustrating Lemma \ref{lem: tsl AxB}.}
\label{fig: AxB}
\end{figure}

\begin{figure}
\centerline{\includegraphics[height=1in]{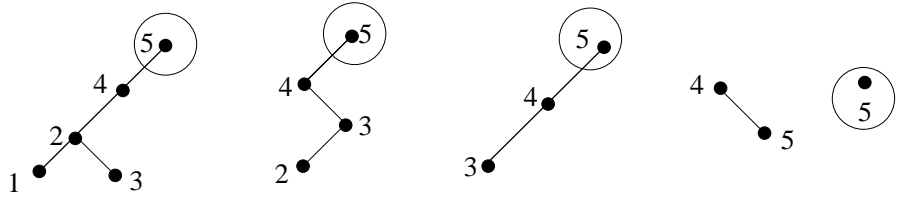}}
\caption{$\rp((M_1)^{[2]}) = 4$. Occurrences of the rightmost
node on the spine are circled.}
\label{fig: m12}
\end{figure}

\begin{definition}
Given a tree $T$ with root $x$, under fetch-and-discard traversal,
the {\em root persistence} of $T$, $\rp(T)$, is the number of steps
$k$, $0 \leq k \leq |T|$, in which $x$ is the root
of $\fetch(k,T)$ (or equivalently, on the spine)
(see Figure \ref{fig: m12}).
\end{definition}

\begin{definition}
The {\em highest echelon} of a tree $T$ is the
rightmost branch leading from the root of $T$.
\end{definition}

\begin{lemma}
\label{lem: tsl AxB}
{\rm (a)} $\tsl(A^xB) = \tsl (A^x) + \tsl(B)$;\hb
{\rm (b)} $\tsl(A^x) = \tsl(A) + \rp(A^x)$;\hb
{\rm (c)} $\tsl(A^xB) = \tsl(A) + \rp (A^x) + \tsl(B)$
(obviously).
\end{lemma}

{\bf Sketch proof.}
(a)
For $0 \leq i \leq |A|$, it can be shown that $\fetch(i,A^x)$ takes
one of the two forms (i,ii) shown in Figure \ref{fig: AxB}.
In either case, $x$ is rightmost on the highest echelon of
$\fetch(i,A^x)$.

In version (i), write $A'$ for $\fetch(i,A)$.  In this case,
\begin{gather*}
A'^x = \fetch(i,A^x),\\
 |\spine(A'^x)| = 1 + \spine(A'), \text{and}\\
A'^xB = \fetch(i,A^xB)
\end{gather*}

In this case, $A^x$ contributes $1$ more unit to $\tsl(A^x)$
than to $\tsl(A)$, and it contributes $1$ more unit
to $\rp(A^x)$.

In version (ii),  there is a node $u$ on the highest
echelon of $\fetch(i,A)$, and in $\fetch(i,A^x)$,
$u$ has right child $x$ and right subtree $X^x$, where
$X$ is the right subtree of $u$ in $\fetch(i,A)$.

In this case, $B$ is the right subtree of $x$ in
$\fetch(i,A^xB$.  Also, $x$ is not at the root,
so this step does not contribute to $\rp(A^x)$.
Also, $\fetch(i,A)$, $\fetch(i,A^x)$, and $\fetch(i,A^xB)$
all have identical spines and the same spine length.

With $i = |A^x|$, $\fetch(i,A^x) = \emptyset$ and
$\fetch(i,A^xB) = B$.  At this point the total
contribution of $A^x$ and $A^xB$ to the total spine
length of both trees is $\tsl(A^x)$.  Continue
the traversal for $|B|$ more steps on $\fetch(|A^x|,A^xB)$
to complete the traversal with total spine length $\tsl(B)$.
Therefore

\be
\tsl(A^xB) = \tsl(A^x) + \tsl(B),
\ee
proving (a).

For (b), having observed that $x$ is on the spine,
and contributes an extra unit to $\tsl(A^x)$ beyond
$\tsl(A)$, and also to $\rp(A^x)$, in case (i) but not case (ii), 
we conclude (b):
\be
\tsl(A^x) = \tsl(A) + \rp(A^x).\qed
\ee

These facts will be used in estimating $\tsl(M_h)$.
Next, the notation $A^x$ will be extended to $A^E$,
where $E$ is a list (ordered) of nodes not in $A$.

\begin{definition}
Let $A$ be a tree and $E = e_1,\ldots, e_k$ a list
of nodes not in $A$.  Inductively one
defines $A^{[E]}$ by: $A^{[\emptyset]} = A$, and for $k>0$
$A^{[e_1,\ldots, e_k]} = (A^{[e_1,\ldots, e_{k-1}]})^{[e_k]}$.

Clearly all trees $A^E$ with $|E|=k$ are isomorphic and
it is often convenient to write $A^{[k]}$ without make the
nodes $e_j$ explicit.
\end{definition}

The following corollary is a version of Lemma
\ref{lem: tsl AxB}, applied to maximal trees.

\begin{corollary}
\be
\label{cor: tsl mh+1 rp mh1}
\tsl(M_{h+1}) = \tsl(M_h^{[1]}) + \tsl(M_h)
=
2 \, \tsl(M_h) + \rp (M_h^{[1]}).\qed
\ee
\end{corollary}

\begin{corollary}
With $\alpha_h$ as defined in Definition \ref{def: alpha second version},
$\alpha_h$ is monotonically increasing, $\alpha$ is well-defined,
and $\alpha_h < \alpha$ for each $h$.
\end{corollary}

{\bf Proof.}
\begin{gather*}
\alpha_{h+1} =
\frac{\tsl(M_{h+1})}{2^{h+2}} = \\
\frac{2 \tsl(M_h) + \rp (M_h^{[1]})}{2^{h+2}} = \\
\frac{\tsl(M_h)}{2^{h+1}} + \frac{\rp(M_{h}^{[1]})}{2^{h+2}} =\\
\alpha_h + \frac{\rp(M_{h}^{[1]})}{2^{h+2}} >
\alpha_h,
\end{gather*}
as claimed. By Elmasry's result the sequence $\alpha_h$ 
is bounded by $8$, so its least upper bound $\alpha$ is well-defined
and for all $h$, $\alpha_h < \alpha$.\qed

\begin{figure}
\centerline{\includegraphics[height=3in]{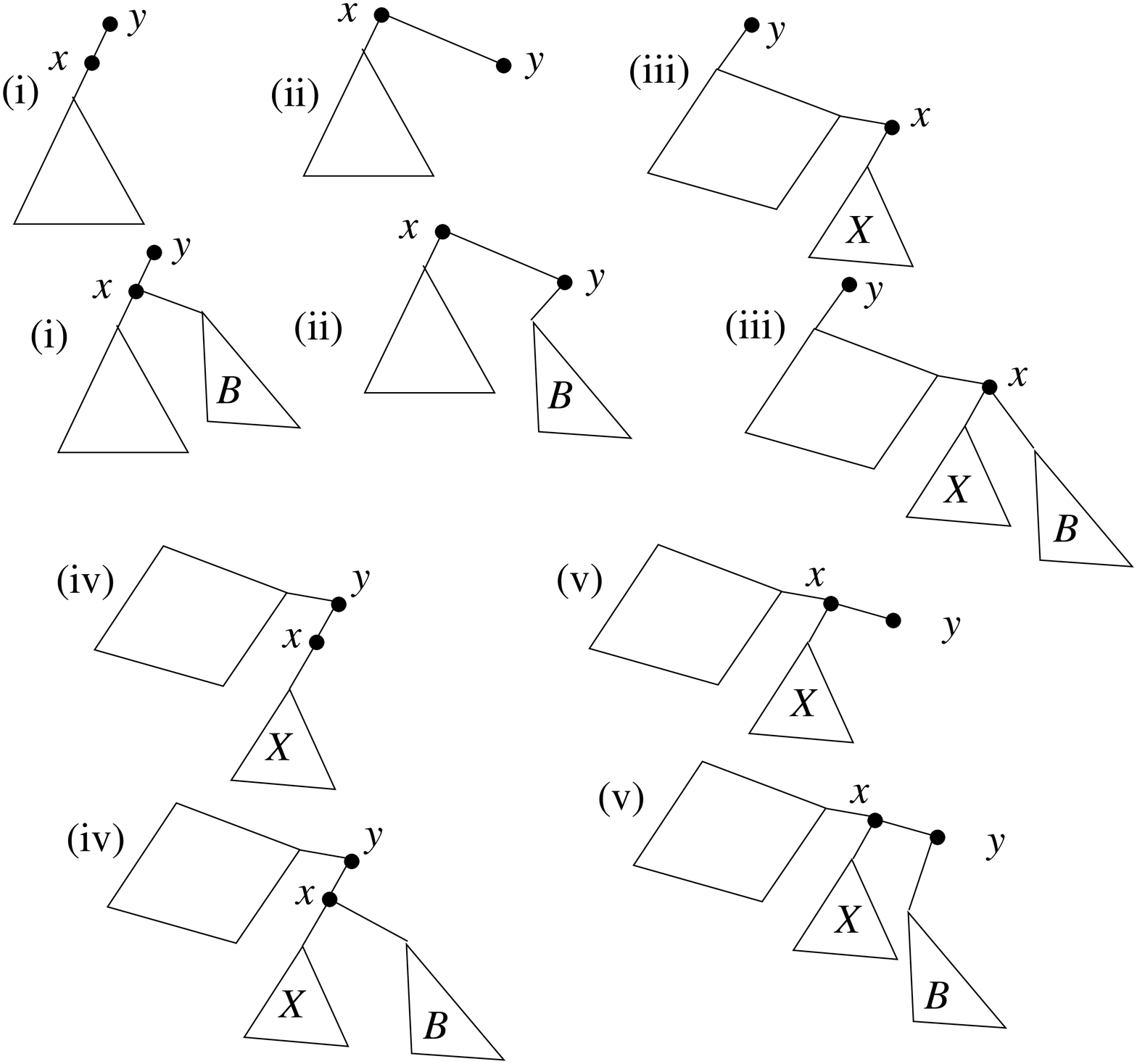}}
\caption{Illustrating Lemma \ref{lem: rp AxBy}.}
\label{fig: AxBy}
\end{figure}

Note that the following lemma is about root persistence, not
total spine length.

\begin{lemma}
\label{lem: rp AxBy}
Given an extended tree $(A^xB)^y$, 
\be
\rp (A^xB)^y = \rp((A^x)^y) - 1 + \rp(B^y).
\ee
\end{lemma}

{\bf Sketch proof.} See Figure \ref{fig: AxBy}.
Also, Figure \ref{fig: m12}.\qed

\begin{corollary}
\label{cor: 2nd diff}
\be
\rp ((M_{h+1})^{[1]}) =
\rp (M_h^{[1]}) + \rp(M_h^{[2]}) - 1.\qed
\ee
\end{corollary}

\section{The initial root persistence of $M_h^{[2]}$}
From Lemma \ref{lem: tsl AxB} (c),
\be
\tsl(M_{h+1}) = 2 \tsl(M_h) + \rp(M_h^{[1]}) .
\ee
So in order to get a fairly sharp upper bound on the
estimate of $\alpha$, we need a fairly sharp upper
bound on $\rp(M_h^{[1]})$, and, in view of the lemma
below, we can use an upper bound  on $\rp(M_h^{[2]})$.

In this section we derive an upper bound on the 
{\em initial} root persistence of $M_h^{[2]}$, which
is the least $k$ such that $y\notin\spine(k,M_h^{[2]})$, where
$y$ is the root and rightmost node of $M_h^{[2]}$
(Lemma \ref{lem: fetching down to 1}).

Since $\rp(M_h^{[2]}) \geq 1$ always,
it follows from
Corollary \ref{cor: 2nd diff}
that $\rp(M_h^{[1]})$ is nondecreasing for $h \geq 0$.
Also, $\rp(M_h)$, since $\rp(M_h^{[1]}) = \rp(M_{h+1})$.

\begin{corollary}
\label{cor: first diff estimate}
\be
\rp(M_h^{[1]}) \leq \sum_{0\leq j < h} \rp(M_j^{[2]}).\qed
\ee
\end{corollary}

{\bf Proof.} Immediate from Corollary \ref{cor: 2nd diff}.
\begin{definition}
If $T$ is a tree and $v$ is the parent of $u$ $\spine(T)$, and a fetch
operation causes $v$ to be made the right child of $u$
(splay operation), we say that $u$ {\em pushes}
$v$ off the spine.

Recall that $\fetch(k,A)$ is the tree after $k$ steps of
fetch-and-discard traversal.

A node $u$ is a {\em repeat node} in $\fetch(k,A)$ if
$u$ is on the spine of that tree, but was pushed off the
spine in a previous step  and later restored to the spine.
\end{definition}

\begin{figure}
\centerline{
\includegraphics[height=2in]{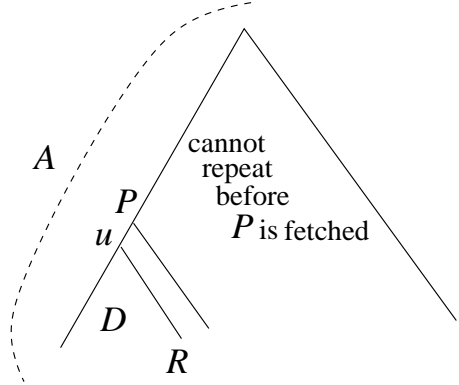}
}
\caption{Illustrating Lemma \ref{lem: no repeat nodes}.}
\label{fig: norepeat}
\end{figure}

\begin{lemma}
\label{lem: no repeat nodes}
Suppose that $u$ is a spine node in $A$ whose
rightmost descendant has inorder rank $k$.
Then there are no repeat  nodes in $\fetch(k,A)$.

Put another way: if $D$ is the subtree with root $u$,
then after  all of $D$ is fetched, there are no repeat nodes.
See Figure \ref{fig: norepeat}.
\end{lemma}

{\bf Proof.}
Let $P$ be the parent of the root of $D$ in $A$, and
$R$ its inorder predecessor, the node in $D$ with
inorder rank $k$.

There were no nodes in the right subtree at $R$ before
$R$ was fetched, for all such nodes would be between $R$
and $P$ in  inorder, and there are none.  Therefore no
nodes are restored when $R$ is fetched.

If $q$ is a repeat node restored before $R$ is fetched,
then it must have been pushed off the spine by its left
child $p$ which is later fetched.   But then $p\in D$ and
therefore $q \in D$, so $q$ is fetched before all of
$D$ is fetched.  Therefore when all of $D$ is fetched,
there are no repeat nodes.\qed

\begin{lemma}
\label{lem: continuation}
Given a tree $A$, let $E$ be an upper segment of $A$
(Definition \ref{def: upper segment}).
Traversing $A$ by fetch-and discard: after
$\lceil \log_2(|E|) \rceil$ fetches, at most $1$ node from 
$E$ has remained continuously on the spine.
\end{lemma}

{\bf Proof.} Let $L_0=E$, and for all relevant $i$ let 
$L_i$ be the set of nodes (in $E$) which have
remained on the spine throughout the first $i$ fetches.
The $(i+1)$-st fetch keeps every second node in $L_i$ on the spine
and pushes the other nodes in $L_i$ off the spine, so by induction,
firstly, $L_{i+1}$ is an interval of contiguous nodes on the spine, 
and secondly

\be
|L_{i+1}|  \leq \left \lceil  \frac{|L_i|}{2} \right \rceil
\ee
so
\be
|L_{i+1}|  \leq \frac{|L_i|}{2} + \frac{1}{2} .
\ee
By induction,
\be
|L_i| \leq \frac{|L_0|}{2^{i}} + 1 - 
\left ( \frac{1}{2}\right )^i .
\ee
Let $i = \lceil \log |E| \rceil$.  Then
\be
|L_i| \leq 1 + 1 - 
\left ( \frac{1}{2}\right )^i < 2
\ee
so $|L_{i}| \leq 1.$\qed

\begin{lemma}
\label{lem: compress spine}
Let $A = M_h^{[E]}$ be an extension of $M_h$.
Let $k = 2 \lceil \log_2 (h+1+|E|) \rceil-1$.
Suppose that $k \leq |M_h|$.  Then within $k$ steps
(or fewer), the spine is reduced to a single node.
\end{lemma}

{\bf Proof.}
Let $k'$ be the smallest number of steps such that
$\spine(\fetch(k',A))$
contains at most one node which has remained continuously
on the spine since the beginning.
From Lemma \ref{lem: continuation}, $ k' \leq \lceil\log_2 (h+1+|E|)\rceil$.
Let $u$ be the smallest subtree of $A$ whose root $u$
is on the spine such that $|D| \geq k'$.  Then $|D| \leq 2k'-1$ (since
$D$ is itself a maximal tree), so $|D| \leq k$,
and when all of $D$ is fetched there is just one node on the spine
(Lemma \ref{lem: no repeat nodes}).\qed

\begin{corollary}
\label{cor: y slips}
In the above lemma, suppose that $y$ is the highest node
of $E$, $|E|\geq 2$, and $k\leq |M_h|$.  Then within $k$
steps, $y$ is pushed off the spine.
\end{corollary}

{\bf Proof.} Let $x\in E$ be the left child of $y$.
Within $k$ steps, the spine is reduced to a single node.
It cannot be $y$, because $y$ is the last node to be fetched,
and no node in $E$ has been fetched.\qed

\begin{definition}
Let $T$ be a tree with root $y$.  The
{\em initial root persistence of $T$}, $\irp(T)$, is
the smallest $k$ such that $y$ is not on the spine
of $\fetch(k,T)$.
\end{definition}

\begin{corollary}
\label{lem: fetching down to 1}
$\irp(M_h^{[2]}) \leq 2 \lceil \log_2 (h+3)\rceil$.\qed
\end{corollary}

\section{Clusters}
The rest of this paper is concerned
with estimating the root persistence of $M_h^{[2]}$.
Recall that the initial root persistence,
$\irp(M_h^{[2]})$, is at most $2\lceil\log_2(h+3)\rceil$.

Because we are interested in the result of fetches near the
top rather than the bottom of trees, the components of
a tree $A^{[2]}$ are labelled as follows.

\noindent (a) The nodes on $\spine(A^{[2]})$ are
labelled $y, x_0, x_1,\ldots$ in {\em descending} order.

\begin{figure}
\centerline{\includegraphics[height=2in]{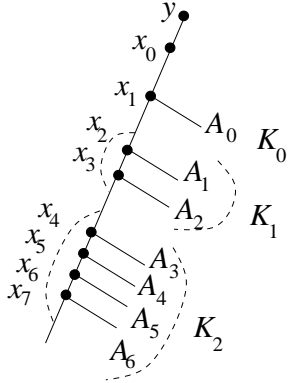}}
\caption{Formation of clusters in $A^{[2]}$.}
\label{fig: clusters}
\end{figure}

\noindent (b)
The right subtrees at $y$ and at $x_0$ are empty.  For any other
spine node $x_i$, $i>0$, the right subtree at $x_i$ is labelled
$A_{i-1}$.
See Figure \ref{fig: clusters}.

The general effect of a single splay operation on a doubly-extended
tree $A^{[2]}$ is as follows:

\begin{itemize}
\item
We call the right branch extending from the left child of
the root the {\em second echelon}.  This is of interest
only when $y$ is on the spine, i.e., $y$ is the root and the only
node on the  highest echelon.

Suppose that the spine is $x_N,x_{N-1},\ldots, x_0,y$ in
{\em bottom to top} order (the indices are decreasing, and $x_N$
is the  lowest vertex on the spine).
\item
$x_N$ is discarded.
\item
$A_N$ is brought onto the spine. That is, if
the spine contains just 1 node ($x_N$) then the tree itself
becomes $A_N$, and otherwise the root of $A_N$ is made
the left child of $x_{N-1}$.
\item
If $N-i$ is odd, and $i\not=0$, then
$x_i,A_{i-1},x_{i-1},A^{i-2}$ are combined into
a single tree with root $x_i$ and right subtree
$(A_{i-1})^{x_{i-1}}A_{i-2}$.

We are interested in the second echelon (when $y$ is on
the spine).  The node $x_0$ is the only node on the
second echelon in $A^{[2]}$, but more generally, as the
traversal proceeds, we may allow more elements on the
second echelon, so $x_0$ has a right subtree $H$, say.

If $N$ is even then $x_1,X_0,x_0,H$ is replaced by
making $(X_0)^{x_0} H$ the right subtree of $x_1$.
This means that $x_1$ joins the second echelon and
$x_0$ is `pushed further along' the second echelon.

\item
If $N$ is odd then $y$ pushed off the spine,
so $H^y$ becomes the right subtree of $x_0$.
\end{itemize}

\begin{figure}
\centerline{\includegraphics[height=2in]{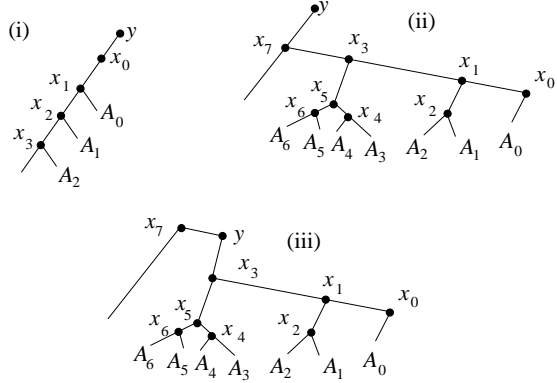}}
\caption{Given $A^{[2]}$ has initial root persistence 4: first,
second, and third clusters are pushed on
the second echelon, then $y$ is pushed off the spine.}
\label{fig: echelon}
\end{figure}

\begin{definition}
\label{def: parts of  fetch A2}
In traversing $A^{[2]}$, suppose that $t < \irp(A^{[2]})$.
Let $A' = \fetch(t,A^{[2]})$: since $t < \irp(A^{[2]})$,
$y$ is the root of $A'$ and has empty right subtree.

The {\em $t$-base} is the smallest subtree $D$ of $A^{[2]}$, {\em not $A'$},
whose root is on the spine of $A^{[2]}$,
and which contains the node of $A^{[2]}$ of inorder
rank $t$.

Looking at $\fetch(t,A^{[2]})$, beginning at the top,
there are
\begin{itemize}
\item
Highest node $y$.
\item
On the second echelon, $t+1$ nodes.  The left subtree
of each $t$ off-spine node is a cluster or partial cluster,
or  possibly neither, a subtree containing nodes only from the base tree $D$.
The leftmost subtree is a $t-1$ cluster or partial cluster
or possibly neither.
\item
Next, some nodes whose right subtree is a $t$-cluster, possibly none.
\item
Next, possibly, a partial cluster.
\item
Then a {\em bottom subtree} containing only nodes from
the $t$-base subtree $D$.
\end{itemize}
\end{definition}

\begin{figure}
\centerline{\includegraphics[height=1in]{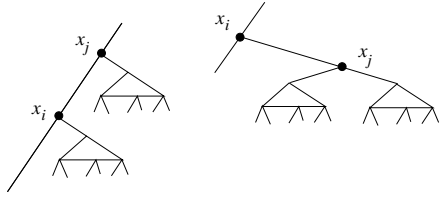}}
\caption{Combining two clusters into one.}
\label{fig: 1-formclusters}
\end{figure}

\begin{figure}
\centerline{\includegraphics[height=2in]{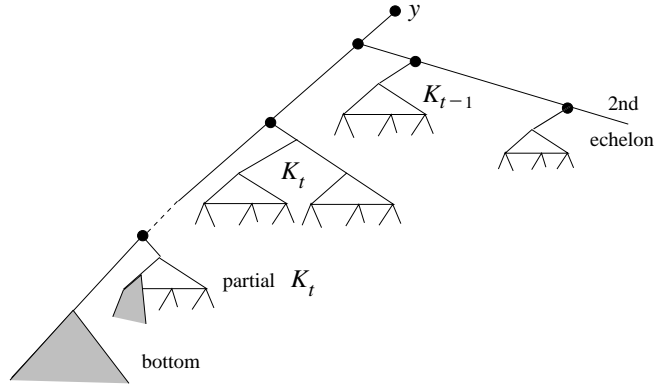}}
\caption{The clusters in $\fetch(t,M_h^{[2]})$, given
$t < \irp(M_h^{[2]}$).}
\label{fig: formechelon}
\end{figure}

\begin{lemma}
If $u,v$ are consecutive nodes on the spine of $\fetch(t-1,M_2^{[2]})$,
where in the next fetch $u$ pushes $v$ off the spine,
{\rm (a)} if the right subtrees of $u$ and $v$ are $t-1$-clusters,
then the combined subtree is a $t$-cluster; {\rm (b)}
if the right subtree of $u$ is a partial cluster,
and $v$ is  not on the second echelon, then the combined
subtree is a partial $t$-cluster; {\rm (c)} if $u$ is  in
the bottom subtree and $v$ is not on
the second echelon and the right subtree of $v$ is a cluster
or partial cluster, then the combined subtree is a partial
cluster; {\rm (d)} if the right subtree of $v$ is from
$D$ then the combined subtree is; {\rm (e)} if $v$ is on
the second echelon, then $u$ joins the second echelon, with
right child $v$, and the right subtree of $u$ becomes the left
subtree of $v$.  In this way another cluster, or partial cluster, or
perhaps part of the bottom subtree, which is before the fetch
right subtree of $u$, becomes the left subtree of $v$.
That subtree is either a $t-1$-cluster, or a partial $t-1$-cluster,
or perhaps neither, being composed entirely of nodes in $D$.

See Figures \ref{fig: 1-formclusters} and
\ref{fig: formechelon}.\qed
\end{lemma}

Below the left and right depth of nodes in a tree are
defined.  This will enable us to prove an important
bound on the length of extensions in traversing a cluster.

\begin{definition}
The {\em depth} of a node $q$ in a tree $T$ is, {\em of course}, the
number of proper ancestors of $q$.  Here we define left and
right depths.

Paradoxically, left depth counts right ancestors and vice-versa.

A {\em right (respectively,
left) ancestor}
of $q$ is a node whose left (respectively, right)
subtree contains $q$.
The {\em left depth} of $q$ is the number of {\em right} ancestors
and the {\em right depth} is the number of {\em left} ancestors.
\end{definition}
Thus the depth is the sum of left and right depths.
If a node is on the  spine, then its left depth is the
number of nodes above it on the spine.

\begin{lemma}
{\rm (i)} Fetch-and-discard on any tree $T$ does not increase the
left depth of any node.  {\rm (ii)} if $T = C^{[2]}$ where
$C$ is a cluster, fetch does not increase the number of
right ancestors which are $x$-nodes.
\end{lemma}

{\bf Proof.}  
(i) Let $q$ be a node before the fetch.  If $q$ is leftmost
then it is fetched and (implicitly) the result follows
automatically. So we assume that $q$ is not leftmost in $T$.

\begin{itemize}
\item
If $q \notin \spine(T)$, say it is in the right subtree
$R$ of a spine node $p$.  Let $S$ be that part of $\spine(T)$
above $p$.

\item
If $q\notin\spine(T)$ and
$p$ is leftmost and fetched, $R$ will be attached to a subsequence
of $S$ which brings $q$ closer to the root of $\fetch(1,T)$.

\item
If $q\notin\spine(T)$ and
$p$ is not leftmost and not pushed off the spine,
then $p$ is brought closer to the root and so is $q$.

\item
If $q\in\spine(T)$ and $q$ is pushed off the spine by
another node $p$, then $q$ acquires a new ancestor $p$,
but it is a left ancestor.

\item
If $q\in\spine(T)$ and is not pushed off the spine, then it
is no further from the root after than before the fetch.
\end{itemize}

(ii): similarly.\qed

\begin{corollary}
\label{cor: |E| leq t+1}
Let $C = K_{t-1}(M_h^{[2]})$ be a cluster.
Suppose that after $f$
fetches, $A_i$ is attached to the spine,
i.e., the root of $A_i$ is on
$\spine(\fetch(f,C^{[2]}))$.  Let $E$ be the
set of spine nodes above this root.  Then
$|E|\leq t+1$.\qed
\end{corollary}

\begin{lemma}
When $A = M_h$,
the $t$-base $D$ (Definition \ref{def: parts of  fetch A2})
has size $|D|\leq 2t-1$. Left subtrees of nodes on the second
echelon which contain only nodes from $D$ have size $\leq |D|$.
In a partial cluster, there are at most $|D|$ nodes
from $D$ (Obviously).\qed
\end{lemma}

\begin{figure}
\centerline{\includegraphics[height=1.5in]{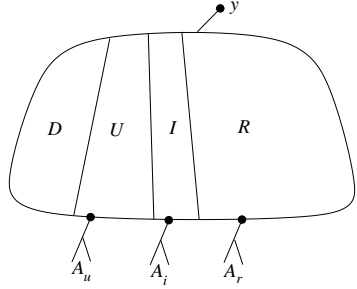}}
\caption{parts of $C^{[2]}$ (see Lemma \ref{lem: rp cluster2}).}
\label{fig: rp-cluster}
\end{figure}

\begin{figure}
\centerline{\includegraphics[height=1in]{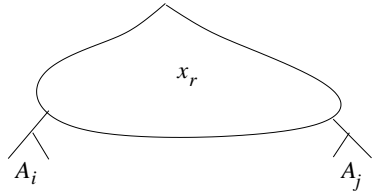}}
\caption{A regular tree $C'$.}
\label{fig: regular-tree}.
\end{figure}

In the lemma below, $t$, and $t \leq \log_2(h+3)$, will be
replaced at one point by $h$, which is assumed to be no smaller.
But $h \geq \log_2(h+3)$ is
only valid for $h\geq 3$.  So, consider $0,1,2$ separately.

\begin{lemma}
\label{lem: h 0 1 2}
\be
\rp(M_0^{[2]} = 2, \rp(M_1^{[2]} = 4,
~\text{\rm and}~ \rp(M_2^{[2]} = 2.
\ee
\end{lemma}

{\bf Sketch proof.}  See Figure \ref{fig: m12} for $\rp(M_1^{[2]})$;
check that $\irp(M_1^{[2]}) = 1$; $\rp(M_0^{[2]}$ is easily checked.

\begin{lemma}
\label{lem: rp cluster2}
If $C$ is a $t$-cluster, or partial cluster,
or pseudo-cluster, derived as part of the
fetch-and-discard traversal of $M_h^{[2]}$,
then
\be
\rp(C^{[2]}) \leq
14 \log_2^2 (4(h+1)) + 2|X| \log_2 (4(h+1))
\ee
where $X$ is the set of $x$-nodes in $C$.
\end{lemma}

{\bf Proof.}  We divide $C$ into four parts,
$D,U,I,R$ as in Figure \ref{fig: rp-cluster}.  Usually
all but the last will be empty.

The  argument is based on the facts that $D$ is small,
$U$ is small ($U$ is based on an index $i_0$ where $h-i_0$
is small), $I$ does not contribute to the root
persistence, and the base trees $A_i$ are sufficiently
large to admit Lemma \ref{lem: compress spine}.

(i) $D$ consists of all nodes in the $t$-base.
An obvious upper bound for the cost of traversing
$D$ is $|D|^2$ which is less than $4t^2$.
Since this counts all steps whether or not $y$ is at
the root, it gives an upper bound for the contribution of these steps
to $\rp(C^{[2]})$.

Let $C' = \fetch(|D|,C^{[2]})$.
Once $D$ is fetched, $C'$ is
{\em regular,} meaning that it consists of an
upper part formed of $x$-nodes, with subtrees $A_i$
attached to it (Figure \ref{fig: regular-tree}).

Let $E$ be the $x$-nodes (plus $y$ if it
is on the spine) on the  spine of $C'$,
and $A_i$ the leftmost $A$-tree, so
$A_i^{[E]}$ is part of $C'$ and $|E| \leq t+1$
(Corollary \ref{cor: |E| leq t+1}).

Now, $A_i = M_{h-i-1}$ with height $h-i-1$.
Let $k = 2 \lceil \log_2(h-i+|E|)\rceil -1$.

From Lemma \ref{lem: compress spine}, if
$k \leq |A_i|$ then within at most
$k$ steps the spine is reduced to a single node,
and either that node is $y$, the only node
left, or it isn't and $y$ has been pushed off the spine
(see Corollary \ref{cor: y slips}).

(ii)
The set $U$ allows for the possibility that
$k > |A_i|$. We focus on those subtrees $A_i$, the $A$-subtrees,
which are attached to 
the regular tree $C'$.

If $k > |A_i|$ then
\begin{gather*}
2 \lceil \log_2(h-i+|E|)\rceil  > 2^{h-i}-1\\
\lceil \log_2(h-i+|E|)\rceil  \geq 2^{h-i-1} \\
 \log_2(h-i+|E|)  > 2^{h-i-1}  - 1
\end{gather*}
The above condition implies the simpler condition
\be
 \log_2(h+t+1)  > 2^{h-i-1} - 1.
\ee
and we {\em define} $U$ by this condition, which
may include more of $C'$ than is necessary, but
leaves nothing out.

If $A_i \in U$ and $i' > i$ then
$A_{i'} \in U$.  So let $i_0$ be
the smallest $i$ satisfying 

\begin{gather*}
 \log_2(h+t+1)  > 2^{h-i-1} - 1\\
2^{h-i-1} < \log_2 ( h+t+1 ) + 1  = \log_2 ( 2(h+t+1) ) \\
h-i-1 <  \log_2 \log_2 ( 2(h+t+1) ) \\
i > h - 1 - \log_2 \log_2 ( 2(h+t+1) ) \\
\text{$i_0$ is minimal:}\\
i_0 = 1 + \lfloor h - 1 - \log_2 \log_2 ( 2(h+t+1) ) \rfloor \\
i_0 = \lfloor h - \log_2 \log_2 ( 2(h+t+1) ) \rfloor
\end{gather*}

The indices $i$ are decreasing and bounded above by $h$.
For each index $i$, there is the tree $A_i$, and the successor
$x$-node,  in $U$.
The range of values for $i \geq i_0$ is at most
\be
\{i: h - i_0 \leq i \leq h\}
\ee

The total number of nodes in $U$ is bounded by
\be
\sum_{i_0}^{h} (1+|A_i|) = \sum 2^{h-i} < 2^{h-i_0+1}
\leq 4 \log_2(2(h+t+1))
\ee
The cost of fetching a node in $U$ is bounded by
$t+1+h-i_0+1$, i.e.,
\be
t+2 + \lceil \log_2\log_2 ( 2 ( h+t+1 ) ) \rceil
\ee

Therefore the total cost of fetching all the nodes in
$U$ is at most

\be
4 \log_2(2(h+t+1)) \times (t+2+\lceil\log_2\log_2(2(h+t+1))\rceil)
\ee

This is an upper bound on the contribution of $U$
to $\rp(C^{[2]})$.

(iii)
The part $I$ allows for $y$ to be restored to
the spine in case after fetching $D$ and $U$
it has been pushed off the spine. There is no
contribution to the root persistence of $C^{[2]}$.

(iv)
Given $C'=\fetch(|D|+|U|+|I|,C^{[2]})$,
$C'$ is regular, and $y$ is on the spine.
Suppose that $A_i$ is the $A$-subtree aligned
with the spine, and $E$ is the remainder of the spine.
Within at most
$k = 2 \lceil \log_2(h-i+|E|)\rceil -1$
fetches, since this time $k$ is not too large,
$y$ is pushed off the spine.

This contributes less than
$ 2 \lceil \log_2(h+t+1)\rceil$
to the root persistence of $C^{[2]}$.

After a certain number of fetches,
which do not contribute to the root persistence,
$y$ is restored to the spine, by a fetch
which removed one of the $x$-nodes.

So this process repeats at most $|X|$ times,
where $X$ is the number of $x$-nodes in $C^{[2]}$,
and the overall root persistence in traversing
$R$ (see Figure \ref{fig: rp-cluster}) is bounded by
\be
2 |X| \lceil \log_2(h+t+1)\rceil
\ee

The total estimate is as follows.
\begin{gather*}
\rp(C^{[2]}) \leq
4t^2 +\\
4 \log_2(2(h+t+1)) \times (t+2+\lceil\log_2\log_2(2(h+t+1))\rceil) +\\
2 |X| \lceil \log_2(h+t+1)\rceil .
\end{gather*}

To simplify this we make some observations.

\begin{gather*}
1 \leq t < \irp(M_h^{[2]}) \leq \lceil \log_2(h+3) \rceil\\
\log_2(h+3) \leq \log_2(2(h+t+2))\\
\text{if}~x\geq 2, \quad \lceil \log_2 x \rceil \leq x\\
\log_2 (2(h+t+1)) \geq 2\\
\lceil \log_2 \log_2 (2(h+t+1)) \rceil \leq
\log_2 (2(h+t+1)) \leq \log_2 (2(h+t+2)).
\end{gather*}
so
\begin{gather*}
\rp(C^{[2]} \leq 4(\log_2^2(2(h+t+2)) +\\
4 \log_2(2(h+t+2)) \times (\log_2(2(h+t+2)) + 2 +
\log_2(2(h+t+2)) +\\
2 |X| \log_2(2(h+t+2)).
\end{gather*}

The root persistence of $M_h^{[2]}$ for
$h = 0,1,2,$ is known
(Lemma \ref{lem: h 0 1 2}): $2,4,2$.  Assuming $h \geq 3$,
$t \leq h$ and we replace $h+t$ by $2h$.

\begin{gather*}
\rp(C^{[2]} \leq 4(\log_2^2(4(h+1)) +\\
4 \log_2(4(h+1)) \times (\log_2(4(h+1)) + 2 +
\log_2(4(h+1)) +\\
2 |X| \log_2(4(h+1)) = \\
12 \log_2^2 (4(h+1)) + 8 \log_2(4(h+1))  +
2|X| \log_2 (4(h+1))
\end{gather*}

\be
\takeanumber
\tag{\thetheorem}
\label{eq: rp cluster2}
\rp(C^{[2]})\leq
14 \log_2^2 (4(h+1)) + 2|X| \log_2 (4(h+1)).
\ee
(This is obviously true when $h\leq 2$.)\qed

Adding these for $1 \leq t \leq 2\log_2 (h+3)$, the
sets $X$ are disjoint and have total size $\leq h+2$, so
\begin{corollary}
\label{cor: first estimate of mh2}
\begin{gather*}
\rp(M_h^{[2]}) \leq
14 \log_2^2 (4(h+1)) (h+3) + 2(h+2) \log_2 (4(h+1))\leq\\
(h+3)(14 \log_2^2 (4(h+1)) + 2\log_2 (4(h+1)) \leq\\
15(h+3) \log_2^2(4(h+1))
\end{gather*}
\end{corollary}

From Corollary \ref{cor: 2nd diff},
\begin{gather*}
\rp(M_h^{[1]} \leq \sum_{j<h} 15(j+3) \log_2^2(4(j+1))
\leq 15 \log_2^2(4(h+1)) (h+3)(h+2)/2\\
\leq 8 (\log_2(h+3)+2)^2 (h+3)^2
\end{gather*}

Now, the estimates of $M_h^{[2]}$ were needed for all
$h \geq 0$, but we
may assume $\log_2(h+3)+2 \leq h+3$ without further investigation
since we need estimates of $M_h^{[1]}$ only for large $h$.
Therefore

\begin{corollary}
\label{cor: estimate of mh1}
For almost all $h$,
\begin{gather*}
\rp(M_h^{[1]}) \leq 8 (h+3)^4.\qed
\end{gather*}
\end{corollary}

\begin{lemma}
\be
\sum_{n\geq N} \frac{n^4}{2^n}  =
\frac{N^4+4N^3+18N^2+52N+75}{2^{N-1}}
\leq
\frac{(N+3)^4}{2^{N-1}}
\ee
\end{lemma}


The equation can be derived easily by assuming that the sum is
a quartic in $N$ divided by $2^N$, and applying the
method of undetermined coefficients.\qed

\begin{corollary}
\label{cor: main inequality}
For any $N>1$,
\be
\alpha \leq \alpha_{N-1} + 8 \frac{(N+3)^4}{2^{N-1}} .\qed
\ee
\end{corollary}

\section*{References}
\begin{enumerate}
\item
\label{elmasry}
Amr  Elmasry (2004).  On the sequential access conjecture
and deque conjecture for splay trees.
{\em Theoretical Computer Science \bf 314:3}, 459--466.


\item
\label{kristensen}
Simon Kristensen (2017).
Arithmetic properties of series of reciprocals of algebraic
integers. Talk delivered to the Department of Mathematics,
Maynooth University, 27 March 2017.




\item
\label{sleator-tarjan}
Daniel Sleator and Robert E. Tarjan (1985).  Self-adjusting
binary search trees. {\em Journal Assoc. Computing Machinery
{\bf 32:3}}, 652--686.



\item
\label{tarjan85}
Robert E.\ Tarjan (1985).
Sequential access in splay trees takes linear time.
{\em Combinatorica \bf 5:4}, 367--378.
\end{enumerate}

\end{document}